\documentclass[12pt,a4paper]{article}
\usepackage{amsmath,amssymb}
\usepackage{bm}

\bmdefine{\aaa}{a}
\bmdefine{\bbb}{b}
\bmdefine{\ccc}{c}
\bmdefine{\eee}{e}
\bmdefine{\ppp}{p}
\bmdefine{\xxx}{x}
\bmdefine{\yyy}{y}
\bmdefine{\zzz}{z}
\bmdefine{\NNN}{N}
\bmdefine{\ZZZ}{Z}
\bmdefine{\QQQ}{Q}
\bmdefine{\RRR}{R}
\bmdefine{\CCC}{C}
\bmdefine{\zerovec}{0}

\newcommand{\aaaa}{{\frak a}}
\newcommand{\bbbb}{{\frak b}}
\newcommand{\mmmm}{{\frak m}}
\newcommand{\nnnn}{{\frak n}}

\makeatletter
\def\ps@outertopright{\let\@mkboth\@gobbletwo%
     \def\@oddhead{\hfill\hbox to 0pt{\hss\thepage\hss}%
     	\hskip -12mm}%
     \let\@evenhead\@oddhead%
     \def\@oddfoot{}\let\@evenfoot\@oddfoot}%
\def\lineup#1{\@ifnextchar[{\@oparglineup{#1}}{\@lineup{#1}}}
\def\@lineup#1#2{\@ifnextchar[{\@@lineupop{#1}{#2}}{\@@lineup{#1}{#2}}}
\def\@@lineup#1#2{\ifmmode%
		#1_1,#1_2,\ldots,#1_{#2}%
	\else%
		$#1_1$, $#1_2$, \ldots,~$#1_{#2}$%
	\fi}
\def\@@lineupop#1#2[#3]{\ifmmode%
		#1_1,#1_2,\ldots,#1_{#2}#3%
	\else%
		$#1_1$, $#1_2$, \ldots,~$#1_{#2}#3$%
	\fi}
\def\@oparglineup#1[#2]#3{\@ifnextchar[{\@@oparglineupop{#1}[#2]{#3}}%
	{\@@oparglineup{#1}[#2]{#3}}}
\def\@@oparglineup#1[#2]#3{\ifmmode%
		#1_{#2},\ldots,#1_{#3}%
	\else
		$#1_{#2}$, \ldots,~$#1_{#3}$%
	\fi}
\def\@@oparglineupop#1[#2]#3[#4]{\ifmmode%
		#1_{#2},\ldots,#1_{#3}#4%
	\else
		$#1_{#2}$, \ldots,~$#1_{#3}#4$%
	\fi}
\def\@@oparglineup#1[#2]#3{\ifmmode\ifinner
			#1_{#2}$, \ldots,~$#1_{#3}%
		\else
			#1_{#2},\ldots,#1_{#3}%
		\fi
	\else
		$#1_{#2}$, \ldots,~$#1_{#3}$%
	\fi}
\def\@@oparglineupop#1[#2]#3[#4]{\ifmmode\ifinner
			#1_{#2}$, \ldots,~$#1_{#3}#4%
		\else
			#1_{#2},\ldots,#1_{#3}#4%
		\fi
	\else
		$#1_{#2}$, \ldots,~$#1_{#3}#4$%
	\fi}
\newcommand{\mysloppy}{\tolerance 9999 \hfuzz .5\p@ \vfuzz .5\p@}
\makeatother

\newtheorem{thm}{Theorem}[section]

\newtheorem{example}[thm]{Example}

\def\define{\mathrel{:=}}

\expandafter\def\csname p@thm\endcsname{(}
\expandafter\def\csname s@thm\endcsname{)}

\newcounter{myequation}[section]

\makeatletter

\def\endequation{\leqno\themyequation$$\global\@ignoretrue}
\makeatother
\def\themyequation{(\arabic{section}.\arabic{myequation})}

\newcommand{\cm}{Cohen-Macaulay}
\newcommand{\rednum}{{\rm r}}
\newcommand{\hilb}{{\rm Hilb}}
\newcommand{\gr}{{\rm Gr}}
\newcommand{\rank}{{\rm rank}}
\newcommand{\height}{{\rm ht}}

\newcommand{\standrep}{standard representation}

\newcommand{\covered}{\mathrel{<\!\!\!\cdot}}
\def\hibiring#1#2{{\cal R}_{#1}(#2)}
\def\int{{\rm int}}

\begin{document}\mysloppy
%\vspace*{2cm}
\begin{center}\LARGE%\huge
On the analytic spread and the
reduction number of the ideal of maximal minors
\end{center}
\begin{center}\Large
Mitsuhiro 
{\sc Miyazaki}%
%\\
%\normalsize
\footnote{Dept. Math. Kyoto University of Education,
\
1 Fukakusa-Fujinomori-cho, Fushimi-ku, Kyoto,
612-8522 Japan,
\
E-mail:\tt
g53448@kyokyo-u.ac.jp}
\end{center}

\begin{trivlist}\item[]\small
{\bf Abstract:}\quad
Let $m$, $n$, \lineup ar, \lineup br\ be
integers with $1\leq a_1<\cdots<a_r\leq m$ 
and $1\leq b_1<\cdots<b_r\leq n$.
And let $x$ be the universal $m\times n$ matrix with
the property that $i$-minors of first $a_i-1$ rows
and first $b_i-1$ columns are all zero, for
$i=1$, \ldots, $r+1$ ($a_{r+1}\define m+1$ and 
$b_{r+1}\define n+1$).
For an integer $u$ with $1\leq u\leq m$,
we denote by $U$ the $u\times n$ matrix consisting of the 
first $u$ rows of $x$.
In this paper, we consider the analytic spread and
the reduction number of the ideal of 
maximal minors of $U$
\end{trivlist}
\begin{trivlist}\item[]\small
{\bf Key words:}\quad
ideal of maximal minors,
analytic spread,
reduction number,
ASL,
distributive lattice
\end{trivlist}

\section{Introduction}
In this paper all rings and algebras are assumed to be
commutative with identity element.
For an $m\times n$ matrix 
$U$ with entries in a ring $R$, we denote by $I_t(U)$
the ideal of $R$ generated by  all the $t$-minors of $U$,
where we put $I_0(U)=R$ and $I_t(U)=(0)$ for $t$ with $t>\min\{m, n\}$.
And if $I_k(U)\neq (0)$ and $I_{k+1}(U)=(0)$, we call
$I_k(U)$ the ideal of maximal minors of $U$.
In this paper, we consider the analytic spread and the reduction
number of the ideal of maximal minors of the matrix
defined below.

Let $K$ be an infinite field, $m$, $n$, $r$ be integers
with $1\leq r\leq \min\{m, n\}$,
\lineup ar, \lineup br\ be integers with $1\leq a_1<\cdots<a_r\leq m$
and
$1\leq b_1<\cdots<b_r\leq n$.
In this situation, there is a universal $m\times n$ matrix
$x$ with the condition
\begin{equation}\label{cond minor vanish}
\left.
\vcenter{%
	\halign{$\displaystyle#$\hfil\cr
	I_i(\mbox{first $a_i-1$ rows})=(0)\cr
	I_i(\mbox{first $b_i-1$ columns})=(0)\cr
	}}\right\}\qquad
\mbox{for $i=1$, \ldots, $r+1$,}
\end{equation}%\label{reduction criterion}
where we set $a_{r+1}\define m+1$ and  $b_{r+1}\define n+1$.
That is, 
$x$ satisfies \ref{cond minor vanish} and
if $U$ is an $m\times n$ matrix with entries in a
$K$-algebra $S$
satisfyings \ref{cond minor vanish}, 
there is a unique $K$-algebra homomorphism
$K[x]\to S$ mapping $x$ to $U$,
where $K[x]$ is the $K$-algebra generated by the entries of $x$.

There are two ways to construct such a matrix $x$.
One is to define $x$ as a homomorphic image of the generic
$m\times n$
matrix (i.e. an $m\times n$ matrix whose entries are independent 
indeterminates) $X$ in the quotient ring of $K[X]$.
The other is to define $x$ as the array of products
of generic matrices, see
\cite{HE}, \cite{HR}, \cite{BV}.
In this paper, we follow the first way.

Let $X$ be the generic $m\times n$ matrix,
$K[X]$ the polynomial ring generated by the
entries of $X$.
Then it is known that $K[X]$ is a graded 
algebra with straightening laws (ASL for short) over $K$
generated by 
$\Delta(X)$,
where
\[
\vcenter{\openup1\jot
	\halign{\hfil$\displaystyle#{}$&$\displaystyle{}#$\hfil\cr
\Delta(X)\define&\{[\lineup cs|\lineup ds]\mid
	s\leq\min\{m,n\},\ 
	\cr
	&\qquad 
	1\leq c_1<\cdots<c_s\leq m,\ 
		1\leq d_1<\cdots<d_s\leq n\}
	\cr
	}
}
\]
and the partial order of $\Delta(X)$ is defined by
\[
\vbox{\halign{&$\displaystyle#$\hfil\cr
[\lineup cs|\lineup ds]\leq[\lineup {c'}{s'}|\lineup {d'}{s'}]\cr
\qquad\Longleftrightarrow
s\geq s',\ %
c_1\leq c'_1,\ldots,c_{s'}\leq c'_{s'},
d_1\leq d'_1,\ldots,d_{s'}\leq d'_{s'}.\cr
}}\]
And $\Delta(X)$ is embedded in $K[X]$ by corresponding
$[\lineup cs|\lineup ds]$ to $\det(X_{c_id_j})_{i,j}$.
See \cite{dep1}, \cite{BV}.

Set
\[
\vcenter{\openup1\jot
	\halign{\hfil$\displaystyle#{}$&$\displaystyle{}#$\hfil\cr
\delta\define&[\lineup ar|\lineup br]\in\Delta(X)\cr
\Delta(X;\delta)\define&\{\gamma\in\Delta(X)\mid \gamma\geq\delta\}\cr
\Omega\define&\Delta(X)\setminus\Delta(X;\delta)\cr
A\define&K[X]/\Omega K[X].\cr
}}
\]
Then, since $\Omega$ is a poset ideal of $\Delta(X)$,
$A$ is a graded ASL over $K$ generated by
$\Delta(X;\delta)$.
If we denote the image of $X$ in $A$ by $x$, then,
by the Laplace expansion, we see that 
$x$ is the universal $m\times n$ matirx
satisfying \ref{cond minor vanish}.

Let $u$ be an integer with $1\leq u\leq m$,
and $U$ be the $u\times n$ matrix consisting of the
first $u$ rows of $x$.
In this paper, we consider the analytic spread
and the reduction number of the ideal of maximal minors $I$
of $U$.
If we take $k$ such that $a_k\leq u<a_{k+1}$,
then by the definition of $A$, we see that
$I_k(U)\neq (0)$ and $I_{k+1}(U)=(0)$.
Therefore $I=I_k(U)$.

%%%%%%%%%%%%%%%%%%%%%%%%%%%%%%%%%%%%%%%%%%%%%%%%%%%%%%%%%%%%%%%%%%%%

\section{Analytic spread}
We denote the irrelevant maximal ideal of $A$ by $\mmmm$.

Northcott-Rees \cite{NR}
defined for an ideal $\aaaa$ of a local ring
$(R,\nnnn)$
with infinite residue field, the 
analytic spread $\ell(\aaaa)$ of $\aaaa$ to be the dimension of
\[
%\ell(\aaaa)\define
R/\nnnn\otimes \gr_\aaaa(R)=
R/\nnnn\oplus \aaaa/\nnnn \aaaa\oplus \aaaa^2/\nnnn \aaaa^2\oplus\cdots
\]
and showed that the analytic spread of $\aaaa$ is the number of
minimal generaters of any
minimal reduction of $\aaaa$.
They also showed, essentially, that for an ideal 
$\bbbb$ contained in $\aaaa$,
\begin{equation}\label{reduction criterion}
\bbbb\aaaa^n=\aaaa^{n+1}\Longleftrightarrow
\overline\bbbb\supseteq(R/\nnnn\otimes\gr_\aaaa(R))_+^{n+1},
\end{equation}%
where $\overline\bbbb$ is the ideal of 
$R/\nnnn\otimes \gr_\aaaa(R)$ generated by
$(\bbbb+\nnnn\aaaa)/\nnnn\aaaa(\subseteq \aaaa/\nnnn\aaaa)$.

Now set
\[
\vcenter{\openup1\jot
	\halign{\hfil$\displaystyle#{}$&$\displaystyle{}#$\hfil\cr
\Theta\define&
	\{\gamma\in\Delta(X)\mid
		[\lineup ak|\lineup bk]\leq\gamma\cr
		&\qquad\qquad\leq[u-k+1,\ldots,u|n-k+1,\ldots,n]\}.\cr
}}
\]

If $\alpha$ and $\beta$ are incomparable elements
in $\Delta(X)$, then the \standrep\ of 
$\alpha\beta$ in the ASL $K[X]$ is of the form
\[
\alpha\beta=\sum_ib_i\gamma_{i1}\gamma_{i2}+\sum_jb'_j\delta_j
\]
and for each $i$ and $j$,
the union of row (column) numbers of 
$\gamma_{i1}$ and $\gamma_{i2}$ ($\delta_j$)
as a multi-set
is the same as that of 
$\alpha$ and $\beta$
\cite{dep1}, \cite{BV}.
Therefore, if $\theta_1$ and $\theta_2$ are incomparable 
elements of $\Theta$, then, since minors of $U$ size 
greater than $k$ are zero,
the standard representation of $\theta_1\theta_2$ in the
ASL $A$ is of the form
\[
\theta_1\theta_2=\sum_ib_i\gamma_{i1}\gamma_{i2},
\qquad\gamma_{ij}\in\Theta.
\]
It follows from \cite[Proposition 1.1]{dep2} that
$K[\Theta]$ is a sub-ASL of $A$
and
\begin{equation}\label{anal sp of A}
A/\mmmm\oplus I/\mmmm I\oplus I^2/\mmmm I^2\oplus\cdots
\simeq K[\Theta].
\end{equation}
Therefore
\[
\ell(I)=\dim K[\Theta]=\rank\Theta+1.
\]
By counting the rank of $\Theta$, we see the following
\begin{thm}
The analytic spread $\ell(I)$ of $I$ is
$k(u+n-k+1)-\sum_{i=1}^k(a_i+b_i)+1$.
\end{thm}

%%%%%%%%%%%%%%%%%%%%%%%%%%%%%%%%%%%%%%%%%%%%%%%%%%%%

\section{Reduction number}
In the following, we multiply the degree of elements in 
$K[\Theta]$ by $1/k$ and adjust the degree of 
the right hand side of \ref{anal sp of A} to the left
hand side.
We also denote the irrelevant maximal ideal of
$K[\Theta]$ by $\nnnn$, and
the analytic spread $\ell(I)$ of $I$ by $l$.
By
\ref{reduction criterion},
if $J$ is  a homogeneous ideal of $A$
generated by elements of degree $k$
(if we consider the degree in $K[\Theta]$, then
degree 1 by the convention above)
and is a minimal reduction of $I$, then the minimal
generating system of $J$ is a homogeneous system of parameters
of degree 1 in $K[\Theta]$.
Conversely, any homogeneous system of parameters of degree 1
in $K[\Theta]$ generates a minimal reduction of $I$ in $A$
(see the proof of \cite[\S2 Theorem 1]{NR}).

If $J$ is a homogeneous ideal of $A$, and a minimal reduction
of $I$, we denote by $\rednum_J(I)$ the reduction 
number $\min\{n\mid JI^n=I^{n+1}\}$ of $I$ with respect to $J$.
If
\lineup vl\ is a minimal system of generators of $J$
of degree 1, then by \ref{reduction criterion},
\[
\vcenter{\openup1\jot
	\halign{\hfil$\displaystyle#{}$&$\displaystyle{}#$\hfil\cr
\rednum_J(I)=&\min\{n\mid(\lineup vl)\supseteq\nnnn^{n+1}\}\cr
	=&\max\{n\mid(K[\Theta]/(\lineup vl))_n\neq0\}\cr
	=&a(K[\Theta]/(\lineup vl)),\cr
}}
\]
where $a(-)$ is the $a$-invariant defined
by Goto-Watanabe (see \cite[Definition (3.1.4)]{GW}).
Since $\Theta$ is a distributive lattice, 
we see that $K[\Theta]$ is a \cm\ ring.
Therefore by
\cite[Remark (3.1.6)]{GW}, we see that
\[
\rednum_J(I)=a(K[\Theta])+l.
\]

On the other hand by the proof of 
\cite[4.4 Theorem]{sta2},
we see that 
\begin{equation}\label{stanley eqn}
\hilb(K_{R},\lambda )=(-1)^{\dim R} \hilb(R,\lambda^{-1})
\end{equation}
for a \cm\ standard graded ring $R$ over a field,
%(note there is a difference of degree in \cite{sta2})
where $\hilb(-,-)$ denotes the Hilbert series,
$K_{R}$ denotes the canonical module of $R$.
So in order to calculate the $a$-invariant, we may replace
the ring with a \cm\ ring with the same Hilbert series.
Since the two ASL's generated by the same poset
has the same Hilbert series, 
we compute the $a$-invariant of $K[\Theta]$ by
computing the $a$-invariant of the Hibi ring
$\hibiring K\Theta$.

In general, for a distributive lattice $D$,
if we denote the set of all the join irreducible elements
of $D$ (i.e.
elements $x$ of $D$ such that there is exactly one $y\in D$
such that $y\covered x$)
by $P$, it is known that
\[
D\simeq J(P)\define\{J\mid\mbox{$J$ is a poset ideal of $P$}\}.
\]
And if one takes a family $\{X_\alpha\}_{\alpha\in P\cup\{-\infty\}}$
of indeterminates and set
$\varphi(I)\define X_{-\infty}\prod_{\alpha\in I}X_\alpha$
for $I\in J(P)$, then Hibi
\cite{hibi}
showed that
\[
{\cal R}_K(D)\define K[\varphi(I)\mid I\in J(P)]
\ (\subseteq K[X_\alpha\mid\alpha\in P\cup\{-\infty\}])
\]
is a homogeneous ASL over $K$ generated by $D$.
Where we set $\deg X_{-\infty}=1$ and $\deg X_\alpha=0$ for
any $\alpha\in P$.

Set
\[
M\define\{(n_\alpha)_{\alpha\in P\cup\{-\infty\}}\in\NNN^{\#P+1}\mid
	\alpha\leq\beta\Longrightarrow n_\alpha\geq n_\beta\}.
\]
Then $M$ is a submonoid of 
$\NNN^{\#P+1}$ and
\[
{\cal R}_K(D)=K[M]\define K[X^\omega\mid \omega\in M],
\]
where $X^\omega$ is the multi-index.

Since $\RRR_{\geq0}M\cap\NNN^{\#P+1}=M$,
we see by
\cite[Theorem 4.1]{sta1},
\[
K_{K[M]}=\bigoplus_{\omega\in\int(\RRR_{\geq0}M)\cap M}KX^\omega.
\]
Because 
$\int(\RRR_{\geq0}M)\cap M=\{
(n_\alpha)_{\alpha\in P\cup\{-\infty\}}\in\NNN^{\#P+1}\mid
	n_\alpha>0$, $
	\alpha<\beta\Longrightarrow n_\alpha> n_\beta\}$,
we see that
\[
\vcenter{\openup1\jot
	\halign{\hfil$\displaystyle#{}$&$\displaystyle{}#$\hfil\cr
a(K[M])=&-\min\{\deg X^\omega\mid\omega\in\int(\RRR_{\geq0}M)\cap M\}\cr
	=&-(\rank P+2).\cr
}}
\]
In particular, by taking $D$ to our $\Theta$, we see
that
\[
\iffalse
\begin{array}{rl}
&\rednum_J(I)\\
=&a(K[\Theta])+l\\
=&a(\hibiring K\Theta)+l\\
=&l-(\rank P+2),
\end{array}
\fi
\rednum_J(I)
=a(K[\Theta])+l
=a(\hibiring K\Theta)+l
=l-(\rank P+2),
\]
where $P$ is the set of join irreducible elements of $\Theta$.

By considering row numbers and column numbers separately,
we see that $\Theta$ is the poset product of two 
distributive lattices, say $D_1$ and $D_2$.
$(x_1,x_2)\in D_1\times D_2$ is join irreducible
if and only if $x_1$ is a join irreducible element of $D_1$
and $x_2$ is the minimal element of $D_2$ or
$x_1$ is the minimal element of $D_1$
and $x_2$ is a join irreducible element of $D_2$.
So if we denote the set of all the join irreducible 
elements of $D_i$ by $P_i$ for $i=1$, $2$, the
set of join irreducible elements of $\Theta$ is 
isomorphic to the disjoint union of $P_1$ and $P_2$.
Therefore
\[
\rank P=\max\{\rank P_1,\rank P_2\}.
\]

Since $D_1$ and $D_2$ are of the same  form,
we consider $D_1$ in the following.
Since
\[
\vcenter{\openup1\jot
	\halign{$\displaystyle#$\hfil\cr
D_1=\{[\lineup ck]\mid
	1\leq c_1<\cdots<c_k\leq u,
	a_i\leq c_i\ (i=1,\ldots,k)\}, \cr
[\lineup ck]\leq[\lineup dk]\Longleftrightarrow\forall i;c_i\leq d_i, \cr
[\lineup ck]\covered[\lineup dk]\Longleftrightarrow
	\exists i;d_i=c_i+1,d_j=c_j(j\neq i),\cr
}}
\]
$[\lineup ck]$ is a join irreducible element of $D_1$
if and only if
there is unique $i$ such that
\begin{equation}\label{join irred}
c_i>a_i,\quad c_i>c_{i-1}+1,
\end{equation}%
where we assume that $c_1>c_0+1$ is always valid.
For a join irreducible element $[\lineup ck]$,
we take $i$ satisfying \ref{join irred},
and set
\[
p\define u-c_i-(k-i),\quad q\define i-1.
\]
We denote the map which send $[\lineup ck]$
to $(p,q)$ by $\varphi$.

It is easy to construct a join irreducible element $[\lineup ck]$
such that $\varphi([\lineup ck])=(p,q)$, if $(p,q)$
is in the image of $\varphi$.
And it also easy to verify that if
$(p,q)$ is in the image of $\varphi$
and $0\leq p'\leq p$, $0\leq q'\leq q$,
then $(p',q')$ is also in the image of $\varphi$.
Moreover,
if
$[\lineup ck]$ and $[\lineup dk]$ are join irreducible elements
of $D_1$,  and
$\varphi([\lineup ck])=(p,q)$, $\varphi([\lineup dk])=(p',q')$,
then
\[
[\lineup ck]\leq[\lineup dk]\Longleftrightarrow
	p\geq p',q\geq q'.
\]
In particular, the coheight of $[\lineup ck]$
in $P_1$ is $p+q$.
Therefore, if we set
\[
%\{i\mid a_i+1<a_{i+1},\ a_i<u,\ i\leq k\}=\{\lineup l[1]v\},\quad
\{i\mid a_i+1<a_{i+1},\ a_i<u,\ i\leq k\}=\{l_1, \ldots, l_v\},\quad
	l_1<\cdots<l_v,
\]
then the minimal elements of $P_1$ are
$[\lineup a{l_1-1},a_{l_1}+1,\lineup a[l_1+1]k]$, \ldots,~%
$[\lineup a{l_v-1},a_{l_v}+1,\lineup a[l_v+1]k]$
and their coheights are
$u-k-a_{l_1}+2l_1-2$, \ldots,~$u-k-a_{l_v}+2l_v-2$
respectively.

\begin{example}\rm
If
$u=13$, $k=8$, $[\lineup ak]=[1,2,3,7,8,10,11,12]$, then
$v=3$, $l_1=3$, $l_2=5$ and $l_3=8$ and
the minimal elements of $P_1$ are
\[
\vcenter{\openup1\jot
	\halign{\hfil$\displaystyle#{}$&$\displaystyle{}#$\hfil\cr
\gamma_1=&[1,2,4,7,8,10,11,12]\cr
\gamma_2=&[1,2,3,7,9,10,11,12]\cr
\gamma_3=&[1,2,3,7,8,10,11,13].\cr
}}
\]
And the Hasse diagram of $P_1$ is the following.
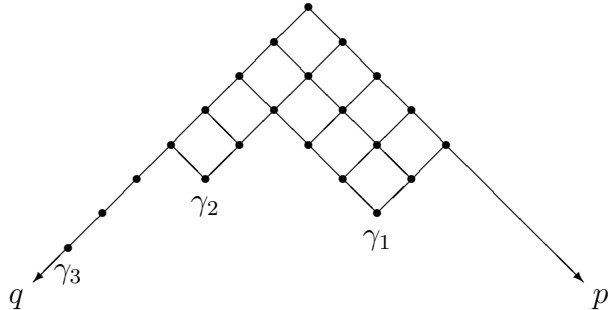
\begin{center}\unitlength\textwidth
\divide\unitlength by 300\relax
\begin{picture}(180,90)(-90,-90)
\put(0,0){\circle*{2}}
\put(-10,-10){\circle*{2}}
\put(-20,-20){\circle*{2}}
\put(-30,-30){\circle*{2}}
\put(-40,-40){\circle*{2}}
\put(-50,-50){\circle*{2}}
\put(-60,-60){\circle*{2}}
\put(-70,-70){\circle*{2}}
\put(10,-10){\circle*{2}}
\put(0,-20){\circle*{2}}
\put(-10,-30){\circle*{2}}
\put(-20,-40){\circle*{2}}
\put(-30,-50){\circle*{2}}
\put(20,-20){\circle*{2}}
\put(10,-30){\circle*{2}}
\put(0,-40){\circle*{2}}
\put(30,-30){\circle*{2}}
\put(20,-40){\circle*{2}}
\put(10,-50){\circle*{2}}
\put(40,-40){\circle*{2}}
\put(30,-50){\circle*{2}}
\put(20,-60){\circle*{2}}

\put(0,0){\vector(-1,-1){80}}
\put(0,0){\vector(1,-1){80}}

\put(10,-10){\line(-1,-1){40}}
\put(20,-20){\line(-1,-1){20}}
\put(30,-30){\line(-1,-1){20}}
\put(40,-40){\line(-1,-1){20}}

\put(-10,-10){\line(1,-1){40}}
\put(-20,-20){\line(1,-1){40}}
\put(-30,-30){\line(1,-1){10}}
\put(-40,-40){\line(1,-1){10}}

\put(-70,-75){\makebox(0,0)[t]{$\gamma_3$}}
\put(-30,-55){\makebox(0,0)[t]{$\gamma_2$}}
\put(20,-65){\makebox(0,0)[t]{$\gamma_1$}}

\put(-85,-85){\makebox(0,0){$q$}}
\put(85,-85){\makebox(0,0){$p$}}
\end{picture}
\end{center}
\end{example}

Summing up, we obtain the following
\begin{thm}
If we set
\[
\vcenter{%
	\halign{$\displaystyle#$\hfil\cr
\{i\mid a_i+1<a_{i+1},\ a_i<u,\ i\leq k\}=\{\lineup l[1]v\},\quad
	l_1<\cdots<l_v
\cr
\{i\mid b_i+1<b_{i+1},\ b_i<n,\ i\leq k\}=\{\lineup {l'}[1]{v'}\},\quad
	l'_1<\cdots<l'_{v'},
\cr}}
\]
then for any minimal reduction $J$ of $I$,
the reduction number $\rednum_J(I)$ of $I$
with respect to $J$ is equal to 
\[
\vcenter{\halign{$#$\hfil&$#$\hfil\cr
\ell(I)-
\max\{&u-k-a_{l_1}+2l_1,\ldots,u-k-a_{l_v}+2l_v,\cr
&n-k-b_{l'_1}+2l'_1,\ldots,n-k-b_{l'_{v'}}+2l'_{v'}\}.\cr}}
\]
\end{thm}


\begin{thebibliography}{DEP2}
%
%
\iffalse
\bibitem[AM]{AM}
Atiyah, M. F. and MacDonald, I. G.:
``Introduction to Commutative Algebra.''
Addison--Wesley (1969)
\fi
%
%
\bibitem[BV]{BV}
Bruns, W. and Vetter, U.:
``Determinantal Rings.''
Lecture Notes in Mathematics {\bf 1327} Springer (1988)
%
%
\bibitem[DEP1]{dep1}
DeConcini, C., Eisenbud, D. and Procesi, C.:
{\it Young Diagrams and Determinantal Varieties.}
Inv.\ Math.\ {\bf56} (1980), 129--165% (1980)
%
%
\bibitem[DEP2]{dep2}
DeConcini, C., Eisenbud, D. and Procesi, C.:
``Hodge Algebras.''
Ast\'{e}risque
{\bf 91} (1982)
%
%
\iffalse
\bibitem[DP]{DP}
DeConcini, C. and Procesi, C.:
{\it A characteristic-free approach to invariant theory,}
Adv.\ Math.\ {\bf21} (1976), 330--354% (1976)
%
%
\bibitem[Fos]{fos}
Fossum, R. M.:
``The divisor class group of a Krull domain.''
Springer (1973)
\fi
%
%
\bibitem[GW]{GW} 
Goto, S. and Watanabe, K.: 
{\it On graded rings, I.}
J. Math.\ Soc.\ Japan {\bf30} (1978), 179--213% (1978)
%
%
\iffalse
\bibitem[GP]{GP}
Gr\"abe, H.-G. and  Pauer, F.:
{\it A remark on Hodge algebras and Gr\"obner bases.}
Czechoslovak Math.\ J. {\bf42} (1992), 331--338% (1992)
\fi
%
%
\iffalse
\bibitem[Hib]{hibi}
Hibi, T.:
{\it Union and glueing of a family of \cm\ partially ordered sets.}
Nagoya Math.\ J. {\bf107} (1987), 91--119% (1987)
\fi
%
%
\bibitem[Hib]{hibi}
Hibi, T.:
{\it Distributive lattices, affine smigroup rings and algebras 
with straightening laws.}
in ``Commutative Algebra and Combinatorics'' (M. Nagata and H. Matsumura, ed.),
Advanced Studies in Pure Math. {\bf11} North-Holland, Amsterdam (1987),
93--109.
%
%
\iffalse
\bibitem[Hoc]{hoc}
Hochster, M.: 
{\it Cohen-Macaulay rings, combinatorics,
and simplicial complexes.}
in ``Ring Theory II,'' Proc.\ of the second 
Oklahoma Conf.\ (B. R. McDonald and R. Morris ed.), Lect.\ Notes in 
Pure and Appl.\ Math., No.26, Dekker (1977), 171--223% (1977)
\fi
%
%
\bibitem[HE]{HE}
Hochster, M. and Eagon, J. A.
{\it \cm\ rings, invariant theory, and the generic perfection 
of determinantal loci.}
Amer. J. Math. {\bf93} (1971), 1020--1058
%
%
\bibitem[HR]{HR}
Hochster, M. and Roberts, J. L.:
{\it Rings of Invariants of Reductive Groups
Acting on Regular Rings are \cm.}
Adv.\ Math.\ {\bf13} (1974), 115--175% (1974)
%
%
\iffalse
\bibitem[Mat]{mat}
Matsumura, H.:
``Commutative ring theory.''
%Cambridge University Press,
%Cambridge studies in advanced mathematics 8 (1986)
Cambridge studies in advanced mathematics {\bf8}
Cambridge University Press
(1986)
\fi
%
%
\iffalse
\bibitem[Miy]{miy}
Miyazaki, M.: 
{\it On 2-Buchsbaum complexes.}
J. Math.\ Kyoto Univ.\ {\bf30} (1990), 367--392% (1990)
\fi
%
%
\iffalse
\bibitem[Mur]{mur}
Murthy, M. P.:
{\it A note on factorial rings.}
Arch.\ Math.\ {\bf15} (1964), 418--420% (1964)
\fi
%
%
\bibitem[NR]{NR}
Northcott, D. G. and Rees, D.:
{\it Reductions of ideals in local rings.}
Proc. Cambridge Phil. Soc. {\bf 50} (1954), 145--158.
%
%
\bibitem[Sta1]{sta1}
Stanley, R. P.:
{\it Linear homogeneous diophantine equations and magic 
labelings of graphs.} Duke Math. J. {\bf40} (1973), 607--632.
%
%
\bibitem[Sta2]{sta2}
Stanley, R. P.:
{\it Hilbert Functions of Graded Algebras.}
Adv. Math. {\bf 28} (1978), 57--83.
%
%
\iffalse
\bibitem[Sta]{sta}
Stanley, R.: 
``Combinatorics and Commutative Algebra.''
Progress in Math.\ {\bf41} Birkh\"auser
%Boston/Basel/ Stuttgart, 
(1983)
\fi
%
%
\end{thebibliography}
\end{document}